\documentclass{amsart}
\usepackage{amssymb,euscript,amsmath, mathrsfs}
\usepackage[dvips]{graphicx}
\usepackage[dvips]{color}
\usepackage{epsfig}

\newcounter{ENUM}
\newcommand{\itm}{\item}
\newenvironment{ilist}{\renewcommand{\theENUM}{\roman{ENUM}}\renewcommand{\itm}{\addtocounter{ENUM}{1}\item[(\theENUM)]}\begin{itemize}\setcounter{ENUM}{0}}{\end{itemize}}

\def\pro{\operatorname{Prop}}
\def\ipro{\operatorname{Improp}}

\def\r{{\bf r}}
\def\S{{\bf S}}
\def\CF{{\mathcal {CF}}}

\def\F{{\mathscr F}}

\def\cH{{\mathcal H}}

\newtheorem{thm}{Theorem}[section]
\newtheorem{prop}[thm]{Proposition}
\newtheorem{lem}[thm]{Lemma}

\theoremstyle{definition}
\newtheorem{defn}[thm]{Definition}

\theoremstyle{remark}

\newtheorem{rem}[thm]{Remark}

\numberwithin{equation}{section}

\begin{document}
\title{Hook length polynomials for plane forests of a certain type}
\author{Fu Liu}
\begin{abstract}
The original motivation for study for hook length polynomials was to
find a combinatorial proof for a hook length formula for binary
trees given by Postnikov, as well as a proof for a hook length
polynomial formula conjectured by Lascoux. In this paper, we define
the hook length polynomial for plane forests of a given degree
sequence type and show it can be factored into a product of linear
forms. Some other enumerative results on forests are also given.
\end{abstract}
\maketitle

\section{Introduction}
In \cite{hook2}, Du and the author defined the hook length
polynomials for $m$-ary trees and showed they can be written as
simple binomial expressions. In this paper, we extend this result to
plane forests of a given degree sequence type.

The original motivation for work on hook length formulas was to seek
a combinatorial proof of an identity derived by Postnikov
\cite{post1, post2}:
\begin{equation}\label{ihookx=1}
\sum_{T}\frac{n!}{2^n}\prod_{v}\left(1+\frac{1}{h_v}\right)=(n+1)^{n-1},
\end{equation}
where the sum is over all complete binary trees with $n$ internal
vertices, the product is over all internal vertices of $T$, and
$h_v$ is the ``hook length" of $v$ in $T$, namely, the number of
internal vertices in the subtree of $T$ rooted at $v$.

Chen and Yang \cite{chenyang} and Seo \cite{seo} both gave direct
bijective proofs for (\ref{ihookx=1}). Moreover, based on
(\ref{ihookx=1}), Lascoux replaced $1$ with $x$ and conjectured a
hook length polynomial formula for binary trees:
\begin{equation}\label{ilas-1}
\sum_{T}\prod_{v}\left(x+\frac{1}{h_v}\right)=\frac{1}{(n+1)!}\prod_{i
= 0}^{n-1} \left((n + 1 + i)x + n +1 -i\right).
\end{equation}

Du and the author \cite{hook2} generalized Lascoux's conjecture and
proved hook length polynomial formulas for $m$-ary trees and plane
forests. Analogous results were also given by Gessel and Seo
\cite{gesselseo}.

In Section 2, we define hook length polynomial a for forests
corresponding to a given degree sequence and show it has a simple
binomial form. In Section 3, we study another form of the hook
length polynomials (\ref{hookp2}) and get an enumerative result on
colored labelled forests (\ref{cfs1}) by using the idea of proper
vertices, which Seo introduced in \cite{seo}. These techniques allow
a fully bijective proof of (\ref{cfs1}), which then yields fully
bijective proofs of the formulas (\ref{hookp}) and (\ref{cpf}) of
section 2.

\section{Hook Length polynomials for plane forests of type $\r$}
A {\it tree} is an acyclic connected graph. For any vertices $v$ and
$u$ in a tree, we call $v$ a {\it descendant} of $u$ (or $u$ an {\it
ancestor} of $v$) if $u$ lies on the unique path from the root to
$v.$ In particular, if $u$ and $v$ are adjacent, we call $v$ a {\it
child} of $u$. For any vertex $v,$ we use $Des(v)$ to denote the set
of descendants of $v.$

 For any vertex $v$ in a tree, the {\it degree} of $v$
is the number of children of $v.$ A vertex is an {\it internal}
vertex if it is not a leaf, i.e., its degree is not zero. We use
$I(F)$ to denote the set of internal vertices of $F.$ A {\it plane
tree} is an unlabelled rooted tree whose vertices are regarded as
indistinguishable, but the subtrees at any vertex are linearly
ordered. A {\it plane forest} is a finite set of ordered plane
trees.

For any plane forest $F,$ let $r_i$ be the number of vertices of
degree $i$ and $\r = (r_0, r_1, r_2, \dots),$ then we say $F$ is of
{\it type} $\r.$ Given a nonnegative integer sequence ${\r} = (r_0,
r_1, r_2, \dots)$ with $\sum_{d \ge 0} r_d < \infty,$ we use
$\F(\r)$ to denote all of the forests $F$ of type $\r.$

There is a well known result about the cardinality of $\F(\r)$
\cite{chen, EE, HPT, stanley2},  denoting by $n = \sum_{d \ge 1} r_d
= |I(F)|$ the number of internal vertices and $\ell = -\sum_{d \ge
0} (d-1)r_d $ the number of trees in $F:$

\begin{equation}\label{cpf}
|\F(\r)| = \frac{\ell}{n + r_0}{{n + r_0} \choose {r_0, r_1, r_2,
\dots}}.
\end{equation}

\begin{defn}For any vertex $v$ of a forest $F,$ we let $d_v$ to be its degree
and $h_v$ its {\it hook length}, i.e. the number of descendants it
has. We define the hook length polynomial of $v$ as
$$P_v(x) = \frac{((d_v - 1)h_v + 1) x + 1 - h_v}{d_v h_v}.$$
\end{defn}

\begin{defn}
We define the hook length polynomial for plane forests of type
$\r$ as:
$$\cH_{\r}(x) = \sum_{F \in \F(\r)} \prod_{v \in I(F)} P_v(x) ,$$

\end{defn}

Then $\cH_{\r}(x)$ can be written as a binomial expression.

\begin{thm}\label{main}
\begin{equation}\label{hookp}
\cH_{\r}(x) = \frac{\ell}{r_0} {{r_0 x} \choose {r_0 x - n, r_1,
r_2, \dots}}.
\end{equation}
\end{thm}

\begin{proof}
If we replace $x$ by $k,$ the right side of (\ref{hookp}) becomes
$$\frac{\ell}{r_0} {{k r_0} \choose {k r_0  - n, r_1, r_2,
\dots}} = \frac{k \ell }{k r_0} {{k r_0} \choose {k r_0 - n + r_1,
r_2, r_3 \dots}}{{k r_0 - n + r_1} \choose r_1}.$$

Applying (\ref{cpf}), one sees that it counts plane forests of type
${\bf r'} = (r_0', r_1', r_2', \dots),$ with $r_1$ leaves circled,
where $r_0' = k r_0 - n + r_1, r_{ik}' = r_{i+1}, \forall i \ge 1 $
and $r_j' = 0$ for all $j \neq ik$ for any $i.$ (Note that $(k r_ 0
- n + r_1) + \sum_{d \ge 2} r_d = k r_0$ and $(k r_0 - n + r_1)
-\sum_{d \ge 2} ((d-1) k - 1) r_{d} = k \ell .$)

Because both sides of (\ref{hookp}) are polynomials in $x,$   it's
enough to prove that
\begin{equation}\label{hook-k}
\cH_{\r}(k) = \frac{\ell}{r_0} {{r_0 k} \choose {r_0 k  - n, r_1,
r_2, \dots}}.
\end{equation}
We prove this by induction on $n,$ the number of internal vertices
of $F.$

When $n = 0,$ we have that $\r = (r_0, 0, 0, \dots)$ and $\ell =
r_0,$ so
$$\cH_{\r}(k) = 1 = \frac{\ell}{r_0}{{r_0 k} \choose {r_0 k}}.$$

Assume (\ref{hook-k}) holds for $n < n_0.$ Now we consider $n =
n_0.$

If $\ell = 1,$ then $\forall F \in \F(\r),$ $F$ is just a tree, say,
$T.$ Let $v_0$ be the root of $T.$ Then
$$\cH_{\r}(k) = \sum_{d \ge 1, r_d \neq 0} \frac{((d-1)n + 1) k + 1 - n}{d n}
\cH_{\r^{(d)}}(k),$$ where $\r^{(d)} = (r_0, r_1, \dots, r_{d-1},
r_d - 1, r_{d+1}, \dots).$

By the induction hypothesis,
\begin{eqnarray*}
\cH_{\r^{(d)}}(k) &=& \frac{d}{r_0} {{r_0 k} \choose {r_0 k  - n +
1, r_1, \dots, r_{d-1}, r_d - 1, r_{d+1}, \dots}} \\
&=& \frac{d}{r_0} \frac{r_d}{r_0 k - n + 1}{{r_0 k} \choose {r_0 k
- n, r_1, \dots, r_{d-1}, r_d, r_{d+1}, \dots}}.
\end{eqnarray*}

Therefore,
\begin{eqnarray*}
\cH_{\r}(k) &=& \sum_{d \ge 1, r_d \neq 0} \frac{((d-1)n + 1) k + 1
- n}{r_0 n} \frac{r_d}{r_0 k - n + 1}{{r_0 k} \choose {r_0 k
- n, r_1, r_2, \dots}} \\
&=& \frac{1}{r_0} {{r_0 k} \choose {r_0 k  - n, r_1, r_2, \dots}}
\sum_{d \ge 1, r_d \neq 0} \frac{(((d-1) n + 1) k + 1 - n) r_d}{ n
(r_0 k - n
+ 1)} \\
&=& \frac{1}{r_0} {{r_0 k} \choose {r_0 k  - n, r_1, r_2, \dots}}.
\end{eqnarray*}

If $\ell > 1,$
\begin{eqnarray*}
\cH_{\r}(k) &=& \sum_{\r = \r^{(1)} + \r^{(2)} + \cdots +
\r^{(\ell)}} \prod_{i=1}^{\ell} \cH_{\r^{(i)}}(k) \\
&=& \sum_{\r = \r^{(1)} + \r^{(2)} + \cdots + \r^{(\ell)}}
\prod_{i=1}^{\ell} \mbox{\# of forests of type ${\r^{(i)}}'$ with
$\r_1^{(i)}$ leaves circled} \\
&=& \mbox{ \# of forests of type $r'$ with $r_1$ leaves
circled} \\
&=& \frac{\ell}{r_0} {{r_0 k} \choose {r_0 k  - n, r_1, r_2,
\dots}}.
\end{eqnarray*}
\end{proof}

\section{Colored Labelled Forests}
In this section, we will use {\it labelled} plane forests. Given a
plane forest $F$ with $n$ internal vertices, a {\it labelling} is a
bijection from $I(F)$ to $[n].$ A {\it labelled forest} is a plane
forest with a labelling. For a vertex $v$ in a labelled forest,
following \cite{seo}, we call $v$ a {\it proper} vertex if none of
its descendants has smaller label than $v,$ and an {\it improper}
vertex otherwise. For a labelled forest $F,$  we use $\pro(F)$ to
denote the set of proper vertices and $\ipro(F)$ the set of improper
vertices.

Suppose we have two sets of colors $\{c_1, c_2, \dots\}$ and $\{
c_1', c_2', \dots \}.$ Fix $k \ge 0,$ a {\it proper $k$-coloring} of
a labelled forest $F$ is a way of color all of the internal vertices
of $F$ so that for any $v \in I(F),$ if $v$ is proper then it can be
colored by any color in $\{c_1, c_2, \dots, c_{d_v} \}$; otherwise
it can be colored by any color in $\{c_1, c_2, \dots, c_{d_v} \}
\cup \{c_1', c_2', \dots, c_k'\}$. (Note that $c_1, c_2, \dots
c_{d_v}$ can be considered corresponding to the $d_v$ edges of $v,$
and the colors $c_1', c_2', \dots, c_k'$ for improper vertices are
considered as ``special'' colors.) Therefore, given a labelled
forest $F,$ there are $\prod_{v \in \pro(F)} d_v \prod_{v \in
\ipro(F)} (d_v + k)$ proper $k$-colorings.

A {\it $k$-colored labelled forest} is a labelled forest with a
proper $k$-coloring. Given a degree sequence $\r$ and $k \ge 0,$ let
$\CF_{\r, k}$ be the set of all $k$-colored labelled forests $F$ of
type $\r.$

\begin{lem}\label{ccf}
$\CF_{\r,k}$ is counted by $\sum_{F \in \F(\r)} n! \prod_{v \in
I(F)}\left((d_v + k) - \frac{k}{h_v}\right).$
\end{lem}

\begin{proof}
For any $F \in \F(\r),$
\begin{eqnarray*}
n! \prod_{v \in I(F)}\left((d_v + k) - \frac{k}{h_v}\right) &=& n!
\sum_{J \subset I(F)} \left( \prod_{v \in J} - \frac{k}{h_v}\right)
\left(\prod_{v \in I(F) \backslash J}(d_v + k) \right) \\
&=&  \sum_{J \subset I(F)} \frac{n!}{\prod_{v \in J} h_v} \left(
\prod_{v \in J} - k\right) \left(\prod_{v \in I(F) \backslash J}(d_v
+ k) \right)
\end{eqnarray*}
However, $n!/\prod_{v \in J} h_v$ is the number of labelling of $F$
so that all the vertices in J are proper. Therefore,
\begin{eqnarray*}
& &n! \prod_{v \in I(F)}\left((d_v + k) - \frac{k}{h_v}\right) \\
&=& \sum_{F'} \sum_{J \subset \pro(F')} \left( \prod_{v \in J} -
k\right) \left(\prod_{v \in I(F) \backslash
J}(d_v + k) \right) \quad (\mbox{$F'$ is $F$ with a labelling}) \\
&=& \sum_{F'} \sum_{J \subset \pro(F')} \left( \prod_{v \in J} -
k\right) \left(\prod_{v \in \pro(F') \backslash J}(d_v + k) \right)
\left(\prod_{v \in \ipro(F')}(d_v +
k) \right) \\
&=& \sum_{F'} \sum_{v \in \pro(F')} \left( -k + d_v + k \right)
\left(\prod_{v \in \ipro(F')}(d_v +
k) \right) \\
&=& \sum_{F'} \left( \mbox{the number of proper $k$-colorings $F'$
has}\right)
\end{eqnarray*}
Summing over all of the forests in $\F(\r)$ gives us:
\begin{equation}\label{clf}
|\CF_{\r,k}| = \sum_{F \in \F(\r)} n! \prod_{v \in I(F)}\left((d_v +
k) - \frac{k}{h_v}\right).
\end{equation}
\end{proof}

Now we look back at our hook length polynomials. We change
(\ref{hookp}) into another form which is more closely related to
Postnikov's identity.

\begin{lem}
 Identity (\ref{hookp}) has the following equivalent form:

\begin{equation}\label{hookp2}
\sum_{F \in \F(\r)} \prod_{v \in I(F)} \frac{(d_v + x)h_v  - x}{d_v
h_v} = \frac{\ell}{r_1! r_2! \cdots} \prod_{i=1}^{n-1}(r_0 + i (1 +
x)).
\end{equation}
\end{lem}

Note when $\r = (n+1, n, 0, 0, \dots)$ and $x = -1,$ (\ref{hookp2})
is the same as (\ref{ihookx=1}).

\begin{proof}
\begin{eqnarray*}
(\ref{hookp}) &\Leftrightarrow& \sum_{F \in \F(\r)} \prod_{v \in
I(F)} \frac{((d_v - 1)h_v + 1) x + 1 - h_v}{d_v h_v} =
\frac{\ell x}{r_1! r_2! \cdots} \prod_{i=1}^{n-1}(r_0 x - i) \\
&\Leftrightarrow& \sum_{F \in \F(\r)} \prod_{v \in I(F)} \frac{((d_v
- 1)h_v + 1) x + (1 - h_v) y}{d_v h_v} = \frac{\ell x}{r_1! r_2!
\cdots}
\prod_{i=1}^{n-1}(r_0 x - i y) \\
&\Leftrightarrow& \sum_{F \in \F(\r)} \prod_{v \in I(F)} \frac{d_v
h_v x + (1 - h_v) (x + y)}{d_v h_v} = \frac{\ell x}{r_1! r_2!
\cdots} \prod_{i=1}^{n-1}(r_0 x - i y) \\
&\Leftrightarrow& \sum_{F \in \F(\r)} \prod_{v \in I(F)} \frac{d_v
h_v y - (1 - h_v) x}{d_v h_v} = \frac{\ell y}{r_1! r_2!
\cdots} \prod_{i=1}^{n-1}(r_0 y + i (x + y)) \\
&\Leftrightarrow& (\ref{hookp2})
\end{eqnarray*}
\end{proof}

If we replace $x$ with $k$ in (\ref{hookp2}) and rearrange it a bit,
we have
\begin{equation}\label{color}
\sum_{F \in \F(\r)} n! \prod_{v \in I(F)}\left((d_v + k) -
\frac{k}{h_v}\right) = {n \choose {r_1, r_2, \dots}}\ell 1^{r_1}
2^{r_2} \cdots \prod_{i=1}^{n-1}(r_0 + i (1 + k)).
\end{equation}

Comparing (\ref{clf}) and (\ref{color}), we get the following
Proposition.

\begin{prop}\label{cf}
$$|\CF_{\r,k}| = {n \choose {r_1, r_2, \dots}}\ell 1^{r_1} 2^{r_2}
\cdots \prod_{i=1}^{n-1}(r_0 + i (1 + k)).$$

\end{prop}

However, we have a stronger result than Proposition \ref{cf}. For
any degree sequence $\r,$ we use $V_{\r}$ to be the set of all
partitions $\S = (S_1, S_2, \dots)$ of $[n]$ such that $|S_i| =
r_i.$ We let $\S = (S_1, S_2, \dots)$ and $\CF_{\r,k,\S}$ be set of
all the forests $F$ in $\CF_{\r, k}$ such that $\forall v \in I(F),$
the label of $v$ is in $S_{d_v}.$

We call two partitions $\S^{(1)}$ and $\S^{(2)}$ {\it adjacent} if
there exists $i \in [n-1],$ such that we can obtain $\S^{(1)}$ by
swapping $i$ and $i+1$ in $\S^{(2)}.$

We construct a graph $G_{\r}$ with vertex set $V_{\r}$ and
$\{\S^{(1)}, \S^{(2)}\}$ forming an edge in $G_{\r}$ if and only if
they are adjacent. It's not hard to see that $G_{\r}$ is connected.

\begin{lem}\label{cfe}
For any two partitions $\S^{(1)}$ and $\S^{(2)}$ in $V_{\r},$
$$|\CF_{\r,k,\S^{(1)}}| = |\CF_{\r,k,\S^{(2)}}|.$$
\end{lem}

\begin{proof}

It's enough to prove the case when $\S^{(1)}$ and $\S^{(2)}$ are
adjacent.

Suppose we obtain $\S^{(1)}$ by swapping $i$ and $i+1$ in
$\S^{(2)},$ for some $i \in [n-1],$ and $i \in S_{d_1}^{(1)}, i+1
\in S_{d_2}^{(1)}.$ (So $i \in S_{d_2}^{(2)}, i+1 \in
S_{d_1}^{(2)}.$)

We define a map $\psi$ from $\CF_{\r,k,\S^{(1)}}$ to
$\CF_{\r,k,\S^{(2)}}.$ For any colored labelled forest $F \in
\CF_{\r,k,\S^{(1)}},$ let $v_1$ be the vertex with label $i$ and
$v_2$ be the vertex with label $i+1$:
\begin{ilist}
\itm If $v_1 \not\in Des(v_2)$ and $v_2 \not\in Des(v_1),$ $\psi(F)
= (i, i+1)
 F,$ where $(i,i+1)F$ means swap the labels $i$ and $i+1.$

\itm If $v_2 \in Des(v_1)$ and $v_1$ is improper, let $\psi(F) = (i,
i+1)F.$

\itm If $v_2 \in Des(v_1)$ and $v_1$ is proper, (so $v_2$ is proper
too,) let $\psi(F) = (i, i+1)F.$

\itm If $v_1 \in Des(v_2)$ and $\exists j < i$ such that the vertex
with label $j$ is in $Des(v_2),$ let $\psi(F) = (i, i+1)F.$

\itm If $v_1 \in Des(v_2),$ any vertex with label $j < i$ is  not in
$Des(v_2),$ and color of $v_2$ is not one of the $k$ special colors,
let $\psi(F) = (i, i+1)F.$

\itm If $v_1 \in Des(v_2),$ any vertex with label $j < i$ is not in
$Des(v_2)$ (so $v_1$ is proper and has a color $c_\alpha$), and the
color of $v_2$ is one of the $k$ special colors, then we obtain
$\psi(F)$ in the following way: Suppose $v_1$ and $v_2$ are in tree
$T$ with root $r.$ Let $u$ be the $\beta$th child of $v_2$ that is
an ancestor of $v_1$ and $w$ be the $\alpha$th child of $v_1$ that
corresponds to the color of $v_1.$ We separate $T$ at $v_2$, $u$,
$v_1$ and $w$ to get five trees $T_1, T_2, T_3, T_4$ and $T_5$ with
roots $r, v_2, u, v_1$ and $w$ respectively, denote by $v_2', u',
v_1'$ and $w'$ the leaves of $T_1, T_2, T_3$ and $T_4$ obtained from
$v_2, u, v_1$ and $w,$ respectively. Attach the root $w$ of $T_5$ to
$u',$ attach the root $v_2$ of $T_2$ to $v_1',$ attach the root $u$
of $T_3$ to $w'$ and attach the root $v_1$ of $T_4$ to $v_2'.$
Color $v_2$ with the color $c_\beta$ corresponding to $u$ and $v_1$
with the original color for $v_2.$ Finally, swap labels $i$ and
$i+1.$

\end{ilist}

One can check that $\psi$ gives a bijection between
$\CF_{\r,k,\S^{(1)}}$ to $\CF_{\r,k,\S^{(2)}}.$
\end{proof}

We observe that Proposition \ref{cf} and Lemma \ref{cfe} together
are equivalent to the following Theorem:

\begin{thm}\label{cfs}
\begin{equation}\label{cfs1}
|\CF_{\r,k,\S}| = \ell 1^{r_1} 2^{r_2} \cdots \prod_{i=1}^{n-1}(r_0
+ i (1 + k)).
\end{equation}
\end{thm}

We also provide another proof of Theorem \ref{cfs}, which is
bijective and combinatorial.

\begin{proof}
Let $\CF_{\r,k,\S,1} \subset \CF_{\r,k,\S}$ be the set with all the
forests with label $1$ appearing in the first tree. Clearly,
(\ref{cfs1}) is equivalent to
\begin{equation}\label{cfs2}
|\CF_{\r,k,\S,1}| = 1^{r_1} 2^{r_2} \cdots \prod_{i=1}^{n-1}(r_0 + i
(1 + k)).
\end{equation}

Let $g_i \in [d_i], \forall 1 \le i \le n, f_j \in [r_0 +i(1+k)],
\forall 1 \le j \le n-1.$ Then there are $1^{r_1} 2^{r_2} \cdots
\prod_{i=1}^{n-1}(r_0 + i (1 + k))$ choices for the $g_i$'s and
$f_j$'s.

We will construct a bijection between $\{g_i, f_j\}$ and
$\CF_{\r,k,\S,1}$ inductively on $n = \sum_{d \ge 1} r_d$ the number
of internal vertices of the forests.

When $n = 1,$ we only have one vertex. Suppose it has degree $d.$
Then $r_0 = d + \ell - 1.$ We don't have $f_j$'s and $g_1 \in [d].$
Clearly, there's a natural bijection between the value of $g_1$ and
the color of vertex $1$ in any forest in $\CF_{\r,k,\S,1}.$

Now we assume for $n < n_0,$ we have a bijection between $\{g_i,
f_j\}$ and $\CF_{\r,k,\S,1},$ and consider $n = n_0.$ For any $F \in
\CF_{\r,k,\S,1},$ let $T$ be it's first tree. We know that $1$ is in
$T.$ We have two cases:

If the root of $T$ is $1,$ then let $g_1$ be the value corresponding
to the color of $1.$ By removing $1$ from $T,$ $F$ becomes a forest
with $\ell + d_1 - 1$ trees and $n - 1$ vertices. However, the
smallest vertex $2$ is not necessarily in the first tree. Let
$f_{n-1}$ be the position number of the tree containing $2$; then
$f_{n-1} \in [\ell + d_1 - 1].$ We cyclicly rotate the order of the
trees so that that tree becomes the first tree in the forest and
call the resulting forest $F'.$ Hence, $F' \in \CF_{\r',k,\S',1},$
where $\r'$ is obtained by subtracting $1$ from $r_1$ in $\r,$ and
$\S'$ is obtained from $\S$ after removing $1$ from $S_{d_1}.$ By
the induction hypothesis, we can associate $\{g_i, f_j\}_{1 \le i
\le n- 1, 1 \le j \le n-2}$ to $F'.$ Including $g_{n}$ and
$f_{n-1}$, we obtain a bijection between forests of this type and
the set $\{g_i, f_j\}$ with $f_{n-1} \in [\ell + d_1 - 1].$

If the root of $T$ is $i \neq 1,$ then $i$ is improper and it can
have $d_i + k$ choices of colors. We can associate these $d_i + k$
colors with choosing $f_{n-1}$ in the interval $[\ell + d_1 +
\sum_{j=2}^{i-1}(d_j+k), \ell + d_1 - 1 + \sum_{j=2}^i (d_j + k)].$
Because $1$ is in $T,$ it is a descendant of $i.$ Let $g_i$ be the
number corresponding to the child of $i$ that is ancestor of $1.$
Similar to the first case, remove $i$ and rotate the first $d_i$
trees so that $1$ becomes contained in the first tree in the new
forest $F'.$ Then $F' \in \CF_{\r',k,\S',1},$ where $\r'$ is
obtained by subtracting $1$ from $r_i$ in $\r,$ and $\S'$ is gotten
from $\S$ after removing $i$ from $S_{d_i}.$ Again, by a similar
argument, we can get a bijection between forests of this type and
the set $\{g_i, f_j\}$ with $f_{n-1} \in [\ell + d_1 +
\sum_{j=2}^{i-1}(d_j+k), \ell + d_1 - 1 + \sum_{j=2}^i (d_j + k)].$

Note that $\ell + d_1 - 1 + \sum_{j=2}^n (d_j + k) = \ell + \sum_{d
\ge 1} (d - 1) r_d + (n - 1)(1 + k) = r_0 + (n-1)(1+k).$ Therefore,
we have constructed a bijection for $n = n_0.$
\end{proof}

\begin{rem}
Because Theorem \ref{cfs} implies Proposition \ref{cf}, Proposition
\ref{cf} together with Lemma \ref{ccf} imply (\ref{hookp2}), and
(\ref{hookp2}) is equivalent to (\ref{hookp}), the above proof of
Theorem \ref{cfs} can be considered as another proof of Theorem
\ref{main}. It also gives a new proof of (\ref{cpf}) by substituting
$-1$ for $x$ in (\ref{hookp}) or substituting $0$ for $x$ in
(\ref{hookp2}).
\end{rem}
\begin{rem}
One can modify our definition of $k$-colorings of a labelled forest
so that it makes sense for $k = -1.$ The proof of Theorem \ref{cfs}
works as well for $k = -1.$
\end{rem}

\bibliographystyle{amsplain}
\bibliography{gen}

\end{document}